\documentclass[runningheads]{llncs}
\usepackage{graphicx}
\usepackage{amsmath,amsfonts}
\usepackage[ruled,linesnumbered,vlined]{algorithm2e}
\usepackage{algorithmic}
\usepackage{stmaryrd}
\usepackage{array}        % Pour personnaliser les colonnes
\usepackage{booktabs, multirow, xcolor, colortbl}
\usepackage{subcaption}  % meilleure alternative à subfigure
\usepackage{url}

\spnewtheorem{assumption}{Assumption}{\bfseries}{\itshape}

\begin{document}
\title{Efficient Solving of Large Single Input Superstate Decomposable Markovian Decision Process}
\titlerunning{Efficient Solving of Large Single Input Superstate Decomposable MDP}

\author{Youssef Ait El Mahjoub \inst{1}\thanks{Corresponding author: youssef.ait-el-mahjoub@efrei.fr} \and
Jean-Michel Fourneau\inst{2,3} \and
Salma Alouah\inst{4} }
\authorrunning{Y. Ait El Mahjoub et al.}
\institute{Efrei Research Lab, Université Paris-Panthéon-Assas, Villejuif, 94800, France
\email{youssef.ait-el-mahjoub@efrei.fr}\\
\and
DAVID Laboratory, Université Paris-Saclay, Versailles, 78000, France
\and 
Inria, ARGO, Paris, France \\
\email{jean-michel.fourneau@uvsq.fr}
\and 
Telecommunications Department, ENSEIRB-MATMECA, Bordeaux INP, France\\
\email{salma.alouah@bordeaux-inp.fr}
}
\maketitle              % typeset the header of the contribution
\begin{abstract}
Solving Markov Decision Processes (MDPs) remains a central challenge in sequential decision-making, especially when dealing with large state spaces and long-term optimization criteria. A key step in Bellman dynamic programming algorithms is the policy evaluation, which becomes computationally demanding in infinite-horizon settings such as average-reward or discounted-reward formulations. In the context of Markov chains, aggregation and disaggregation techniques have for a long time been used to reduce complexity by exploiting structural decompositions. In this work, we extend these principles to a structured class of MDPs. We define the Single-Input Superstate Decomposable Markov Decision Process (SISDMDP), which combines Chiu’s single-input decomposition with Robertazzi’s single-cycle recurrence property. When a policy induces this structure, the resulting transition graph can be decomposed into interacting components with centralized recurrence. We develop an exact and efficient policy evaluation method based on this structure. This yields a scalable solution applicable to both average and discounted reward MDPs.

\end{abstract}
%Please note that the first paragraph of a section or subsection is
%not indented. The first paragraph that follows a table, figure,
%equation etc. does not need an indent, either.
%Subsequent paragraphs, however, are indented.

\section{Introduction}
Solving large-scale Markov chains remains a fundamental challenge in a variety of domains, including performance evaluation of computer systems, reliability analysis, and biological modeling. As the state space grows, classical exact methods become computationally prohibitive, particularly when attempting to compute stationary probability distributions or long-term performance metrics. This has led to the development of a wide range of methods around aggregation and disaggregation \cite{courtois1977,Stew94,buchholz1994}, which aim to aggregate the original Markov chain into a smaller system that can be solved more efficiently, followed by a refinement or reconstruction phase. These methods leverage structural properties such as lumpability \cite{buchholz1994} and quasy-lumpability \cite{Muntz1994,Marin2022} or weakly connected components (NCD - Near Completely Decomposable Markov Chains), and have become standard tools for analyzing Markovian systems.

When extending this setting to Markov Decision Processes (MDPs), the computational cost grows considerably, as each action introduces its own transition model, thereby increasing the overall complexity of decision-making. However, even for a fixed policy, where the MDP reduces to a single Markov chain, evaluation remains computationally expensive, particularly in infinite-horizon formulations such as average-reward or discounted-reward criteria. In such cases, policy evaluation typically involves solving large linear systems or performing iterative updates, and must be repeated multiple times within dynamic programming algorithms (e.g., policy iteration, value iteration).

To overcome this, several structured MDP frameworks have been proposed to exploit regularities in the model. Hierarchical MDPs (HMDPs) \cite{barto2003,dietterich2000} decompose decision problems into nested sub-tasks or options \cite{sutton1999}, enabling abstraction and reuse of sub-policies in a temporally extended decision process. In contrast, Factored MDPs (FMDPs) \cite{boutilier1995,koller1999} focus on compact representations of the state space, modeling it as a product of variables and leveraging conditional independence to factor transition and reward models, which enables efficient inference and policy computation in high-dimensional domains. While these structured frameworks focus respectively on functional decomposition in the case of HMDPs and probabilistic factorization in the case of FMDPs, our approach introduces a fundamentally different form of structure based on the topology of the transition graph induced by the policy. Rather than constraining the action space or explicitly defining variable dependencies, we exploit a structural organization that emerges naturally from the dynamics of the policy. This organization combines localized recurrence with constrained inter-component communication, giving rise to a new class of decision processes with internal regularities that can be exploited for efficient computation.

The structural foundations of our approach build upon two classical models from the theory of Markov chains. Chiu et al. \cite{Chiu87} introduced the notion of Single-Input Superstate Decomposable Markov Chains (SISDMC), in which the state space is partitioned into strongly connected components, and all transitions between components must enter through a unique designated root state. In parallel, Robertazzi \cite{Rob90} studied chains where all internal cycles are constrained to pass through a central root, enforcing a form of centralized recurrence within each component. We have previously demonstrated the effectiveness of the Robertazzi model in both purely stochastic and decision-based contexts. Specifically, in \cite{YHJM18,YHJM19}, we employed this structure to model the filling process of an optical container. In \cite{YAJM24}, we extended it to a Markov Decision Process for modeling the energy filling process in a battery station under stationary energy arrivals. This was further generalized in \cite{YAJM25}, where we considered non-stationary arrivals driven by photovoltaic (PV) panel production, requiring a more dynamic control policy. As a natural continuation of \cite{YAJM25}, we turned to Chiu's structure, which can be viewed as a generalization of Robertazzi’s model by allowing inter-component communication through root states, making it particularly suitable for modeling multi-station systems. However, this paper does not primarily address the multi-station context, but instead introduces an efficient policy evaluation algorithm for a class of MDPs that integrate both structural properties. We define the resulting model as a \textit{Single-Input Superstate Decomposable Markov Decision Process (SISDMDP)}, where each partition satisfies Robertazzi’s single-cycle condition, and the global interconnection follows Chiu’s single-input topology. The main contribution of this paper is to show that such structured MDPs admit a fast and exact policy evaluation method, grounded in the recursive decomposition of their transition graph.

The remainder of the paper is organized as follows. Section~\ref{sec:model} introduces the SISDMC-SC structure in the context of Markov chains and presents the proposed SISDMDP model. Section~\ref{sec:resolution} details the resolution of these models under both the average and discounted reward criteria, along with complexity analysis. Section~\ref{sec:results} provides numerical results comparing the proposed method to standard algorithms. Finally, Section~\ref{sec:conclusion} concludes the paper and discusses future directions.
\vspace{-0.5cm}
\section{Model Description \label{sec:model}}
\vspace{-0.3cm}
We begin by defining the SISDMC-SC (Single Input Superstate Decomposable Markov Chain – Single Cycle) structure. Consider an irreductible Markov chain with $N$ states.

\begin{definition} \cite{Chiu87}
A \textit{single-input superstate} is a subset of states in a Markov chain such that exactly one state within the subset receives incoming transitions from states outside the subset. This state is called the \emph{input state} (or superstate). All other states in the subset can only be reached from within the subset itself. Formally, let $\mathcal{S} = \{1, 2, \dots, N\}$ be the state space of a Markov chain with transition matrix $P = (P_{i,j})$. A subset $S = \{s_1, s_2, \dots, s_m\} \subset \mathcal{S}$ is called a single-input superstate with input state $s_1$ if
\vspace{-0.15cm}
\[
P_{j,i} = 0 \quad \forall j \notin S, \; i \in S \setminus \{s_1\}. \vspace{-0.2cm}
\]
\end{definition}

\begin{definition} \cite{Chiu87}
A \textit{Single-Input Superstate Decomposable Markov Chain (SISDMC)} is a Markov chain that can be divided into multiple disjoint superstates, each of which satisfies the single-input condition. Formally, the state space $\mathcal{S}$ can be partitioned as $\mathcal{S} = \bigsqcup_{i=1}^K S_i$ where $K \geq 2$, and each $S_i$ is a single-input superstate.
\end{definition}

\begin{definition} \cite{Rob90}
In a \textit{Rob-B structure} (as defined by Robertazzi), every directed cycle in the Markov chain passes through a single specific state. %This reflects a centralized control or bottleneck structure within the dynamics.
\end{definition}

\begin{definition}
We define a new structure called the \textit{SISDMC-SC} model, which combines the SISDMC structure with the Rob-b cycle constraint. In this model, each partition must not only satisfy the single-input property but also enforce that all internal cycles go through the superstate state of that partition.
\end{definition}

\begin{lemma}
\label{lemaGen}
    By definition, the SISDMC-SC structure is a generalization of the Rob-B structure, specifically when considering $K = 1 $ partition.
\end{lemma}

\begin{figure}[ht]
\hspace{-1cm}
  \begin{subfigure}[b]{0.6\textwidth}
    \includegraphics[height=5.5cm,width=\textwidth]{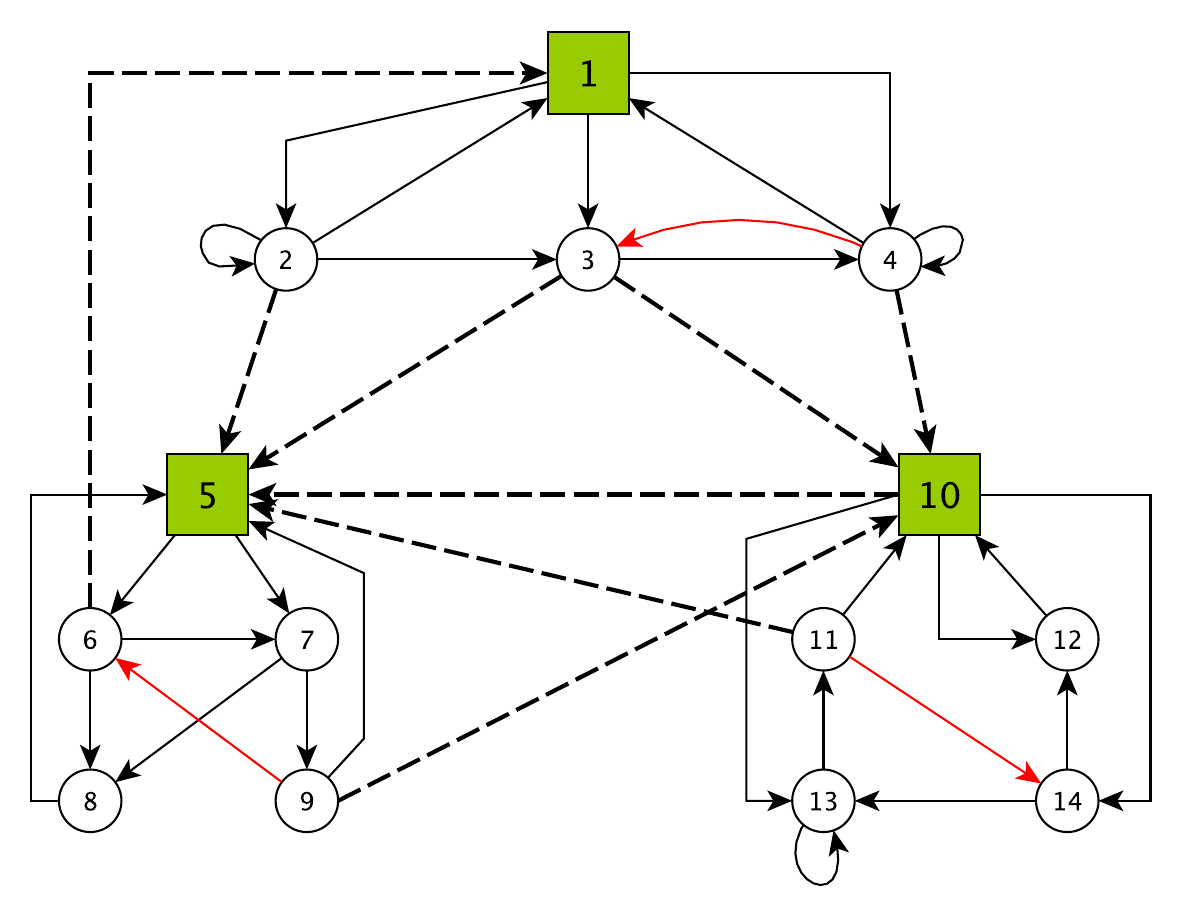}
    \vspace{-0.4cm}
    \caption{SIDMC stucture}
    \label{fig:partitionA}
  \end{subfigure}
  \hspace{0.5cm}
  \begin{subfigure}[b]{0.6\textwidth}
\includegraphics[height=5.5cm,width=\textwidth]{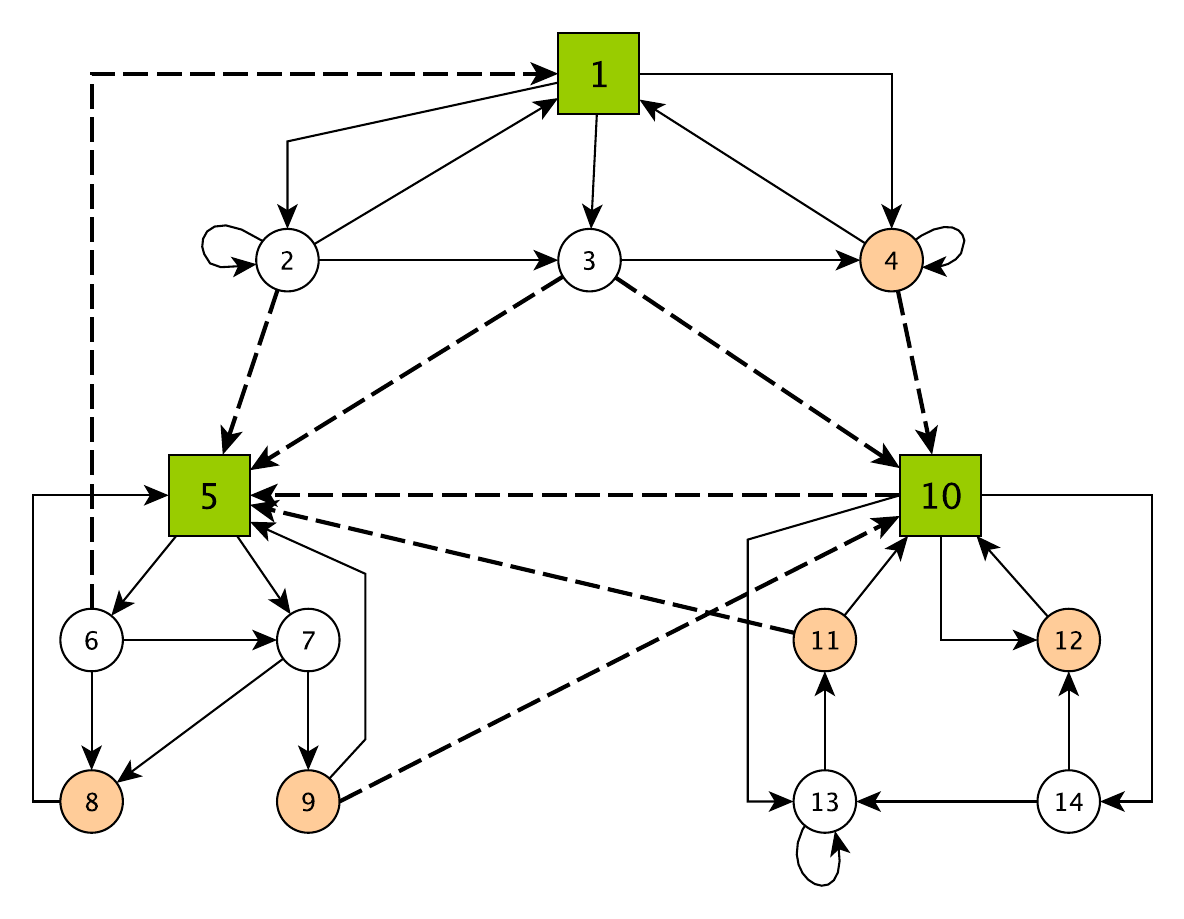}
    \vspace{-0.4cm}
    \caption{SISDMC-SC structure}
    \label{fig:partitionB}
  \end{subfigure}
  \label{structure}
\end{figure}

In Fig.~\ref{fig:partitionA}, we illustrate an example of a SISDMC model. Green states correspond to superstates. Solid arcs represent transitions within the same partition, while dotted arcs denote inter-partition transitions. Note that in Fig.~\ref{fig:partitionB}, the SISDMC-SC structure is obtained from the original SISDMC by removing the red arcs, as they can form cycles that do not pass through a superstate.

This work generalizes the use of this structural pattern from Markov chain analysis to the resolution of decision-making problems within the framework of Markov Decision Processes (MDPs). An MDP is defined as a tuple $\{S,\ A,\ P^{(a)},\ R^{(a)}\}$, where $S$ is a finite set of states, $A$ a finite set of actions, $P^{(a)}$ the transition probability matrix for action $a$, and $R^{(a)}$ the immediate reward function associated with transitions under action $a$.

A \emph{policy} $\pi : S \to A$ is a mapping from states to actions, specifying the action to be taken in each state. The objective is to determine an optimal policy $\pi^*$ that maximizes the expected reward over an infinite decision-making horizon. In particular, we focus on two standard formulations: maximizing the \emph{discounted cumulative reward} and the \emph{long-run average reward}.

We now formally define the class of MDPs considered in this work:

\begin{definition}
\label{SISDMDP}
A \emph{Single-Input Superstate-Decomposable Markov Decision Process (SISDMDP)} is an MDP such that, for any policy $\pi$, the transition graph induced by $P^{(\pi)}$ exhibits a SISDMC-SC structure. Additionally, we assume that the resulting Markov chain is ergodic.
\end{definition}
\vspace{-0.6cm}
\section{MDP resolution \label{sec:resolution}}
\vspace{-0.3cm}
In following, we leverage the structural pattern of the SISDMDP to evaluate efficiently, with exact results, any policy in the  evaluation phase of the policy iteration algorithm \cite{Putr94,Abhj15}. 
\vspace{-0.3cm}
\subsection{Average reward criteria}
\vspace{-0.2cm}
First, we recall Bellman equations in the context of Policy Evaluation algorithm.  To optimize an average reward criteria, under policy $\pi$, one need to estimate
\begin{itemize}
    \item the average reward $\rho^{(\pi)} $, representing the expected reward per time step,
    \item and the relative value function $V^{(\pi)} $, capturing deviations from this average in each state.
\end{itemize}
The average reward starting from state s is defined as 
\vspace{-0.3cm}
\begin{equation}
    \rho^{(\pi)}(s) = \lim_{T \to \infty} \frac{1}{T} \sum_{t=0}^{T-1} \mathbb{E}^{(\pi)} \left[ r(s_t, \pi(s_t)) \,\middle|\, s_0 = s \right], \vspace{-0.1cm}
\end{equation}
where  $r(s_t,\pi(s_t))$ is the immediate reward obtained from state $s$ at time $t$ taking action $\pi(s_t)$.
In practice, a simplification occurs when the Markov chain induced by policy $\pi $ is \emph{unichain} \cite{Putr94}. That is the average reward does not depend on the state. In unchain policies, the induced graph generates a single recurrent class (with some transient states). Hence, states will be revisited indefinitely which leads, asymptotically, to similar average reward. Unlike multi-chain policies, which can generate multiple recurrence classes, resulting in a possible distinct average value for each recurrence class.
\begin{lemma}
\label{lemUnichain}
    The SISDMDP is unichain: \vspace{-0.15cm}
    \begin{equation}
    \forall s, s', \ \ \rho^{(\pi)}(s) = \rho^{(\pi)}(s') = \rho^{(\pi)}. \vspace{-0.15cm}
\end{equation}
\end{lemma}
\begin{proof}
By assumption, for every stationary policy 
$\pi$, the induced SISDMC-SC is ergodic, that is, irreducible and aperiodic. Consequently, the Markov chain induced by any such policy contains a single recurrent class that includes all states. This implies that the SISDMDP is unichain. Therefore, the average reward obtained from a decision trajectory starting from any initial state converges to the same value, denoted $\rho^{(\pi)}$.
\end{proof}

Next, we introduce the value function $V^{(\pi)}$ associated with policy $\pi$, defined as the cumulative expected reward starting from state $s$. However, the natural value function in average reward criteria tends to diverge unless we subtract the average reward. This contrasts with the discounted reward, where discount factor $\gamma<1$ ensures to have bounded values from estimated future rewards. 
A natural version of the value function is defined as, $\forall s\in S$ as  
\vspace{-0.3cm}
\begin{equation}
     V^{(\pi)}(s) = \lim_{T \to \infty} \mathbb{E}^{(\pi)} \big[ \sum_{t=0}^{T-1} r(s_t,\pi(s_t)) \ | \ s_0 = s \big]. \vspace{-0.15cm}
\end{equation} 
This expression is equivalent, in matrix form, to $V^{(\pi)} (I - P^{(\pi)\top}) = R^{(\pi)}$ where $P^{(\pi)\top} $ is the transpose of the transition matrix $P^{(\pi)} $, and $R^{(\pi)} $ is the reward vector under policy $\pi $.  However, this system is difficult to solve directly because $I - P^{(\pi)\top} $ is a singular matrix. This singularity reflects the divergence of values often encountered in the average reward framework. To address this problem, a relative value function is defined which consists in retrieving the value function of some defined state $x$ (i.e. the relative value), solving the singularity issue. Hence $\forall s\in S$
\vspace{-0.4cm}
\begin{equation}
     V^{(\pi)}(s) = \lim_{T \to \infty} \mathbb{E}^{(\pi)} \big[ \sum_{t=0}^{T-1} r(s_t,\pi(s_t)) \ | \ s_0 = s \big] -  V^{(\pi)}(x) \vspace{-0.4cm}
\end{equation}
\begin{equation}
\label{eqV}
      \Rightarrow \ \ V^{(\pi)}(s) = r(s,\pi(s)) -  \rho^{(\pi)} +  \sum_{s'=1}^{N}P^{(\pi)}_{s,s'} \ V^{(\pi)}(s').
\end{equation}

This last formulation is the Bellman equation for relative  policy evaluation \cite{Putr94,Abhj15} which consists on a system of $N$ linear equations. The unknowns are vector $V^{(\pi)}$ and  scalar $\rho^{(\pi)}$. That could be either solved by classical linear solvers that comes with significant computational cost  or iteratively, with some lack of precision, using fixed point methods. One note that if the steady-state distribution for some policy, we note $\Pi^{(\pi)}$, exists then we can derive the average Markov reward process formula \vspace{-0.3cm}
\begin{equation}
\label{eqRho}
    \rho^{(\pi)} = \sum_{s\in \mathbb{S}} \Pi^{(\pi)}(s).r(s,\pi(s)) . \vspace{-0.1cm}
\end{equation}
Once a policy is evaluated (i.e. by solving equation system \eqref{eqV}), one can use following equations to improve the policy. The Q-function is defined as \vspace{-0.2cm}
\begin{equation}
\label{eqQ}
Q(s,a) = r(s,a) + \lambda \sum_{s'=1}^{N}P^{(a)}_{s,s'}V(s'),
\vspace{-0.1cm}
\end{equation}
hence optimal policy \cite{Putr94,Abhj15} in each state is defined as \vspace{-0.2cm}
\begin{equation}
\label{eqpi*}
    \pi^*(s) \in \arg\max_{a\in A(s)}\Big[ Q(s,a) \Big].\vspace{-0.1cm}
\end{equation}
(Note that $\lambda=1$ for the average reward criteria)

The Relative Policy Iteration (RPI) algorithm begins with an arbitrary policy, which is evaluated using Equation~\eqref{eqV}. The policy is then improved, if possible, using Equations~\eqref{eqQ} and~\eqref{eqpi*}. The algorithm stops when no further improvement is possible according to Equation~\eqref{eqpi*}; the resulting policy is then the optimal policy $\pi^*$.

In this work, our goal is to solve Equation~\eqref{eqV} efficiently for the SISDMDP class, as this represents the most time-consuming phase of the RPI algorithm. To that end, we must first compute the average reward $\rho^{(\pi)}$, which requires obtaining the steady-state distribution $\Pi^{(\pi)}$. Once $\rho^{(\pi)}$ is calculated, it can be substituted into Equation~\eqref{eqV} to complete the policy evaluation step, which we also aim to accelerate by exploiting the structural properties of the SISDMDP.

We first recall that in the Rob-B topology, there are two main types of intra-superstate structures, that is states can be ordrer such that: \vspace{-0.1cm}
\begin{itemize}
    \item $P^{(\pi)} = C^{(\pi)} + U^{(\pi)}$, where $C^{(\pi)}$ is a matrix whose first column is positive and all other entries are zero, and $U^{(\pi)} $ is an upper triangular matrix. The first state, $s_1$, corresponds to the root of the subgraph, that is, the superstate of the intra-superstate structure. The resulting graph is an arborescence with return cycles directed back to the superstate (typically, partitions $\{1,\dots,4\}$ and $\{5,\dots,9\}$ in Fig.~\ref{fig:partitionB}). This type of structure can model filling or accumulation processes, such as those observed in optical containers~\cite{YHJM18} or battery charging dynamics~\cite{YAJM24,YAJM25}.
    
    \item $P^{(\pi)} = D^{(\pi)} + L^{(\pi)} $, where $D^{(\pi)} $ is a matrix whose first row is positive and all other entries are zero, and $L^{(\pi)} $ is a lower triangular matrix. The resulting graph is then an anti-arborescence (typically, partition $\{10,\dots,14\}$ in Fig.~\ref{fig:partitionB}). This structure is suited to representing data collection networks, such as LoRa-based sensor systems \cite{SLTD17,LAZM15}. \vspace{-0.3cm}
\end{itemize}

\begin{remark}[Partition types]
In the remainder of this paper, we assume that partitions follow the first structure. This assumption is made for clarity and without loss of generality: our analysis and techniques readily extend to the second case, and more generally, to any SISDMDP instance involving a mixture of both types of partitions.
\end{remark} 
\vspace{-0.4cm}

\subsubsection{I- Calculating $\rho^{(\pi)}$:} In \cite{Chiu87}, Chiu and Feinberg presented an efficient and direct algorithm for computing the steady-state probability distribution of SISDMC Markov chains. The key idea is to isolate each partition of the state space by redirecting external transitions to the superstate of the corresponding partition. This results in what is known as the \emph{intra-superstate system}. The steady-state distribution is first computed locally within each partition. Then, a reduced \emph{inter-superstate system} is constructed by considering only the superstates and their interactions. Finally, the global steady-state distribution is obtained by combining the local (intra-superstate) and global (inter-superstate) steady-state vectors via a vector product.
However, Chiu's method does not specify which numerical algorithm should be used to solve each subsystem (e.g., GTH \cite{GTHe85}, Power method, Gauss-Jordan elimination, etc.). In our case, for the SISDMC-SC structure, we take advantage of the Rob-B topology and apply the efficient algorithm, Algorithm \ref{algo-Rob}, which we have previously validated in \cite{YHJM18,YHJM19}. This algorithm solves each intra-superstate system in linear time, with complexity $O(m)$, where $m$ denotes the number of arcs (non-zero transitions). The adapted version of Chiu's method for the SISDMC-SC structure is presented in Algorithm \ref{algo-ChiuNew}.

\begin{algorithm}[ht]
\SetAlgoLined
\SetKwInOut{Input}{Input}
\SetKwInOut{Output}{Output}
 \caption{\label{algo-Rob}Steady-state algorithm for policy $\pi$, Rob-B model}
 \Input{Transition matrix $P^{(\pi)}$, superstate $s_1$}
 \Output{Steady-state probability distribution vector $\Pi^{(\pi)}$}
 \BlankLine
     Initialize $\alpha(s_1) = 1$   \\
     Get values of $\alpha(s)$ for all $s>s_1$ using  
     \ \ \  \ \ \  \ \ \ $\alpha (s)  = \big(\sum_{s'<s} \alpha(s') U^{(\pi)}[s',s]\big)/ \big(1-U[s,s]),  \ \ \forall s>s_1 $ \\
     Deduce the value of $\Pi^{(\pi)}(s_1)$ from normalisation  $\Pi^{(\pi)} (s_1)  = \big[ 1+ \sum_{s>s_1}^{N}\alpha (s) \Big]^{-1} $  \\
     Obtain $\Pi^{(\pi)}(s)$ for all $s>s_1$ using $\Pi^{(\pi)} (s)  = \alpha (s) \ \Pi^{(\pi)} (s_1)$.
\end{algorithm}

For clarity in the presentation of Algorithm~\ref{algo-ChiuNew}, we introduce the following notation:  
Let $K $ be the number of partitions. Each partition $S_r$ (with $r \in \llbracket 1, K \rrbracket $) contains $n_r $ states, and we denote by $S_r = \{s_{1,r}, \dots, s_{n_r,r}\} $ the set of states in partition $S_r$. The first state $s_{1,r} $ represents the superstate of $S_r$. %Let $\mathit{S_{sup}} $ denote the set of all superstates, i.e.,  
%\[
%\mathit{S_{sup}} = \{s_{1,r} \mid r \in \llbracket 1, K \rrbracket\}.
%\]

\begin{algorithm}[ht]
\SetAlgoLined
\SetKwInOut{Input}{Input}
\SetKwInOut{Output}{Output}
\caption{\label{algo-ChiuNew} Average reward for policy $\pi$, SISDMC-SC model}
\Input{Transition matrix $P^{(\pi)}$, instant reward $r(s,\pi(s))$, partitions $\{S_r\}_{r=1}^K $}
\Output{Average reward $\rho^{(\pi)}$ induced by policy $\pi$}
\BlankLine

\ForEach{partition $S_r = \{s_{1,r}, \dots, s_{n_r,r}\} $}{
  - Construct the intra-superstate matrix $A^{(\pi)}_r \in \mathbb{R}^{n_r \times n_r} $ as
\begin{equation}
      \forall i,j \in \llbracket 1, n_r\rrbracket, \ \ \  A^{(\pi)}_r(i,j) =
  \begin{cases}
    P^{(\pi)}(s_{i,r}, s_{j,r}) & \text{if } s_{j,r} \ne s_{1,r} \\
    1 - \sum\limits_{k=2}^{n_r} P^{(\pi)}(s_{i,r}, s_{k,r}) & \text{if } s_{j,r} = s_{1,r}
  \end{cases}
\end{equation}
  -Execute \textbf{\textit{Algorithm \ref{algo-Rob}}} to obtain $\phi^{(\pi)}_r$ the steady-state probability vector induced by $A^{(\pi)}_r$
}
Construct the inter-superstate matrix $B^{(\pi)} \in \mathbb{R}^{K \times K} $ as 
\begin{equation}
B^{(\pi)}(r,q) = 
\begin{cases}
    \sum\limits_{i=1}^{n_r} \phi_r^{(\pi)}(i) \cdot P^{(\pi)}(s_{i,r}, s_{1,q}) & \text{for} \ r \neq q \\
    1 - \sum\limits_{q \neq r} B^{(\pi)}(r,q) & \text{for} \ r = q
\end{cases}
\end{equation}
This represents the transition probability from the superstate $s_{1,r} $ to the superstate $s_{1,q} $ of partition $q $, weighted by the steady-state vector $\phi_r^{(\pi)} $. \\

Obtain $\psi^{(\pi)}$ the steady-state probability vector induced by $B^{(\pi)}$ \\

Derive the overall steady-state probability of $P^{(\pi)}$
\begin{equation}
\Pi^{(\pi)} = \left[ \psi_1^{(\pi)} \cdot \phi_1^{(\pi)}, \psi_2^{(\pi)} \cdot \phi_2^{(\pi)}, \dots, \psi_K^{(\pi)} \cdot \phi_K^{(\pi)} \right]
\end{equation}\\
Deduce the average reward $\rho^{(\pi)}$ from Equation \eqref{eqRho} using $r(s,\pi(s))$ and $\Pi^{(\pi)}$. 
\end{algorithm}
\vspace{-0.3cm}
\begin{lemma}
\label{lemComplexite1}
The complexity of Algorithm~\ref{algo-ChiuNew} is \vspace{-0.2cm}
\begin{equation}
\mathcal{O}\left( \sum_{r=1}^{K} m_r + K^3 \right),
\end{equation}
 The first term accounts for the local steady-state computations using Algorithm~\ref{algo-Rob}, which runs in linear time with respect to the number of arcs. The second term corresponds to the resolution of the global steady-state system over the inter-superstate matrix $B^{(\pi)} $, which is dense and solved using the GTH algorithm \cite{GTHe85} with cubic complexity.

This approach is significantly more efficient than the classical Chiu method, which applies a cubic-cost solver such as GTH to each local matrix $A^{(\pi)}_r $, resulting in a total cost of
$\mathcal{O}\left( \sum_{r=1}^{K} n_r^3 + K^3 \right).$
\end{lemma}
\vspace{-0.4cm}
\subsubsection{II- Calculating $V^{(\pi)}$:}
To compute $V^{(\pi)} $, it is important to distinguish between calculating the steady-state system and the value-state system. In the former, we compute probabilities, and more specifically, in the so-called balance equations, each state is expressed as a function of \textit{incoming transitions} (i.e., $\Pi^{(\pi)} P^{(\pi)} = \Pi^{(\pi)} $). In contrast, in the value-state system (Equation \eqref{eqV}), which is part of a decision-making formulation involving real values rather than probabilities, the equivalent matrix equation is $V^{(\pi)} = V^{(\pi)} \cdot P^{(\pi)\top} + R^{(\pi)} $, where each state is expressed as a function of \textit{outgoing transitions}, which are derived from the transposed matrix of $P^{(\pi)} $.

Linear systems can be solved using classical direct methods with cubic complexity, or via fixed-point iterative approaches with $\mathcal{O}(iter.N^2)$ complexity (\textit{iter} is the number of iterations needed for convergence), the latter may suffer from limited numerical precision. In \cite{YAJM24,YAJM25}, we proposed an efficient method for solving the value-state system in the context of the Rob-B structure to model a battery filling process with intermittent energy arrivals. The method relies on the property of the unique "bias" state (or relative state) in relative value evaluation. This property ensures that the estimated value of the relative state is fixed at 0. To solve the system, we fix the bias state as the root state and perform a bottom-up propagation throughout the system to efficiently deduce the values of other states. This method is efficient in structures with a common entry point or isolated partitions. However, it does not extend to the SISDMC-SC structure, where multiple partitions communicate with each other. This discrepancy leads us to consider an alternative approach: keeping the same reasoning as in the former model, that is, if we estimate the values of all superstates, we can propagate their values within each partition to obtain the values of all other states. 

It is important to note that each state can either transition to intra-partition states, with a unique possible cycle passing through the superstate, or transition directly to other superstates (by definition). This implies that each state can be expressed as a function of all the superstates. Hence, the first step of our method is to derive the inter-superstates system (i.e., a linear system composed solely of superstates). By solving this system, we can propagate the values within each partition. Let’s define the following sets. $S_{\text{sup}} $ as the set of all superstates: \vspace{-0.15cm}
\begin{equation}
    S_{\text{sup}} = \{ s_{1,1}, s_{1,2}, \dots, s_{1,K} \}.
    \vspace{-0.1cm}
\end{equation}
Next, we define $S_{r,R} $ (resp. $S_{r, \overline{R}} $), the set of Release (resp. non-Release) states in partition $S_r $ that have transitions only to superstates (resp. states that can transition to both superstates and other states within the partition), as follows:\vspace{-0.15cm}
\begin{equation}
\label{eqSets}
    \begin{cases}
        S_{r,R} = \left\{ s_i \in S_r \mid \forall s_j \in S_{r, \ j \neq 1}\ ,\ P^{(\pi)}(s_i, s_j) = 0 \text{ and } \exists\ s_k \in S_{\text{sup}},\ , P^{(\pi)}(s_i, s_k) > 0 \right\} \\
        S_{r, \overline{R}} = S_r \setminus S_{r,R}
    \end{cases}
\end{equation}
In Fig.~\ref{fig:partitionB}, we illustrate the sets $S_{r,R}$ in light red. Specifically, $S_{1,R} = \{4\}$, $S_{2,R} = \{8, 9\}$, and $S_{3,R} = \{11, 12\}$. These states are essential, as they mark the starting point of the substitution procedure.
\vspace{-0.3cm}
\subsubsection{II-A) Local Substitution Within Partitions:}
Let us now describe the construction of the linear system involving only the superstates. The key idea is to \emph{eliminate the non-superstates by expressing their values as linear combinations of the values of superstates}, exploiting the structure of the model. For each partition $r \in \llbracket 1, K \rrbracket $, we derive a system of the form:
\begin{equation}
    M^{(r)} V_{\text{sup}} = b^{(r)},
\end{equation}
where $V_{\text{sup}} \in \mathbb{R}^K$ is the vector of unknown values for all superstates $S_{\text{sup}}$ , and the matrix $ M^{(r)} \in \mathbb{R}^{n_r \times K}$, and the vector  $b^{(r)} \in \mathbb{R}^{n_r}$, are built recursively. We now construct $M^{(r)}$  and $ b^{(r)}$  in two phases based on the internal structure of the partition.

\paragraph{1. Release States ($S_{r,R}$):}

For all $s_i \in S_{r,R}$, the state has no transitions to other intra-partition states (except possibly to the superstate). Hence,  Equation \eqref{eqV} simplifies to one involving only superstates: \vspace{-0.15cm}
\begin{equation}
\label{eqVNR}
 V^{(\pi)}(s_i) = \frac{1}{d(s_i)} \left( r(s_i,\pi(s_i)) -  \rho^{(\pi)} +  \sum_{s_j \in S_{\text{sup}}} P^{(\pi)}(s_i, s_j) V(s_j) \right) 
\end{equation}
where $d(s_i) = 1 - P^{(\pi)}(s_i, s_i)$  that corresponds to possible self-loops in $s_i$ state.
Following Equation \eqref{eqVNR}. In this step, the $i$-th row of $M^{(r)}$ denoted as $M^{(r)}_{i,\cdot}$ will store the normalized transition probabilities toward superstates, and $b^{(r)}_i$ stores the normalized immediate reward minus average reward:
\begin{equation}
    \label{eqM1}
        M^{(r)}_{i,\cdot} \gets \frac{1}{d(s_i)} \cdot \left( P^{(\pi)}(s_i, s_j) \right)_{s_j \in S_{\text{sup}}}
    \end{equation}
\begin{equation}
\label{eqB1}
    b^{(r)}_i \gets \frac{r(s_i,\pi(s_i)) - \rho^{(\pi)} }{d(s_i)}. \vspace{-0.3cm}
\end{equation}

\paragraph{2. Non-Release States ($S_{r, \overline{R}}$):}

For these states, the Bellman equation includes contributions from both intra-partition transitions and transitions to superstates. To handle these, we \emph{recursively substitute} the equations of previously treated states, with a bottom-up procedure as $S_{r, R} $ are bottom states in a partition. Let $s_i \in S_{r, \overline{R}} $, its value can be written as: \vspace{-0.3cm}
\begin{equation}
\label{eqVR}
    V(s_i) = \frac{1}{d(s_i)} \left( r(s_i,\pi(s_i)) -  \rho^{(\pi)} + \sum_{s_j \in S_r \setminus \{s_i\}} P^{(\pi)}(s_i, s_j) V(s_j) + \sum_{s_k \in S_{\text{sup}}} P^{(\pi)}(s_i, s_k) V(s_k) \right).
\end{equation}

The term $V(s_j) $ for $s_j \in S_r \setminus \{s_i\} $ is substituted using the rows already built in $M^{(r)} $ and $b^{(r)} $. The full bottom-up substitution induced by Equation \eqref{eqVR} gives: \vspace{-0.3cm}
\begin{equation}
    \label{eqM2}
    M^{(r)}_{i,\cdot} \gets \frac{1}{d(s_i)} \left( \sum_{s_j \in S_r \setminus \{s_i\}} P^{(\pi)}(s_i, s_j) M^{(r)}_{j,\cdot} + \sum_{s_k \in S_{\text{sup}}} P^{(\pi)}(s_i, s_k) e_k \right),
\end{equation}
\begin{equation}
\label{eqB2}
    b^{(r)}_i \gets \frac{1}{d(s_i)} \left( r(s_i,\pi(s_i)) -  \rho^{(\pi)}  + \sum_{s_j \in S_r \setminus \{s_i\}} P^{(\pi)}(s_i, s_j) b^{(r)}_j \right),
\end{equation}
where $e_k \in \mathbb{R}^K $ denotes the canonical basis vector with a 1 in the $k $-th coordinate (corresponding to the superstate $s_k $) and 0 elsewhere. Transitions to superstates are thus handled exclusively in $M^{(r)} $, ensuring that each state is expressed as a linear function of superstate values only. \vspace{-0.4cm}
\subsubsection*{II-B) Global System Extraction:}
From each local system $M^{(r)} V_{\text{sup}} = b^{(r)} $, we extract the equation corresponding to the root (i.e., the superstate $s_{1,r} $), which is always located in the first row of $M^{(r)} $. This yields a global system involving only the superstates: \vspace{-0.1cm}
\begin{equation}
\label{eq:Vsup1}
A V_{\text{sup}} = B,
\vspace{-0.1cm}
\end{equation}
where
\begin{equation}
\label{eq:Vsup2}
A = \begin{bmatrix}
M^{(1)}_{1,\cdot} \\
M^{(2)}_{1,\cdot} \\
\vdots \\
M^{(K)}_{1,\cdot}
\end{bmatrix}
\in \mathbb{R}^{K \times K},
\quad
B = \begin{bmatrix}
b^{(1)}_1 \\
b^{(2)}_1 \\
\vdots \\
b^{(K)}_1
\end{bmatrix}
\in \mathbb{R}^{K}.
\vspace{-0.45cm}
\end{equation}
\subsubsection*{II-C) Resolution:}
To remain consistent with the relative policy evaluation, we fix the value of a reference superstate (e.g., $V(s_{1,1}) = 0 $). The resulting linear system can then be solved using any classical method (e.g., Gauss-Jordan elimination), yielding the vector $V_{\text{sup}} $ of superstate values to be propagated back into each partition. \vspace{-0.3cm}
\subsubsection*{II-D) Final Injection:}
Once the values of the superstates $V_{\text{sup}} $ are known, we propagate them within each partition to reconstruct the full value function $V^{(\pi)} \in \mathbb{R}^N $. The value of each superstate $s_{1,r} $ is already known from the solution of the superstates system and is directly injected into $V^{(\pi)} $. For the remaining states $s_{i,r} \in S_r \setminus \{s_{1,r}\} $, their values are reconstructed using the local system $M^{(r)} V_{\text{sup}} + b^{(r)} $, as follows: \vspace{-0.4cm}
\begin{equation}
\label{eqReconstruct}
V(s_{i,r}) = \sum_{k=1}^{K} M^{(r)}_{i,k} \cdot V_{\text{sup}}[k] + b^{(r)}_i \vspace{-0.2cm}
\end{equation}
This step completes the policy evaluation under the relative value formulation. 
\vspace{-0.4cm}
\subsection{Discounted reward criteria}
\vspace{-0.2cm}
In contrast to the average reward setting, the discounted reward formulation focuses on maximizing the cumulative reward obtained over time, while discounting future rewards with a factor $\gamma \in [0, 1[ $. Under a fixed policy $\pi $, the value function associated with the discounted criterion is defined as: \vspace{-0.3cm}
\begin{equation}
V^{(\pi)}(s) = \mathbb{E}^{(\pi)} \big[ \sum_{t=0}^{\infty} \gamma^t r(s_t, \pi(s_t)) \ | \  s_0 = s \big]. \vspace{-0.3cm}
\end{equation}

Unlike the average reward case, the natural discounted formulation does not require ergodicity or unichain assumptions. The existence of the value function $V^{(\pi)} $ are guaranteed as long as the reward function is bounded and $\gamma < 1 $. This makes the discounted criterion particularly appealing for theoretical analysis and for algorithms relying on contraction properties. The value function $V^{(\pi)} $ also satisfies the Bellman fixed-point equation: $\forall s \in S$ \vspace{-0.15cm}
\begin{equation}
\label{eqVdisc}
V^{(\pi)}(s) = r(s, \pi(s)) + \gamma \sum_{s' \in \mathcal{S}} P^{(\pi)}_{s,s'} \, V^{(\pi)}(s'). \vspace{-0.3cm}
\end{equation}
The proposed policy evaluation procedure remains structurally identical to that used in the average reward setting. However, several adjustments are required to account for the discounted formulation. First, it is no longer necessary to compute the average reward $\rho^{(\pi)} $. Second, when evaluating $V^{(\pi)} $, the transition matrix must be scaled as $P^{(\pi)} \gets \gamma \cdot P^{(\pi)} $; that is, each entry of $P^{(\pi)} $ is multiplied by the discount factor $\gamma $ (this substitution is applied in step A).

In addition, the definition of the vector $b^{(r)}_i $ in Equation~\eqref{eqB1} must be updated as follows: \vspace{-0.25cm}
\begin{equation}
\label{eqdB1}
    b^{(r)}_i \gets \frac{r(s_i, \pi(s_i))}{d(s_i)}, \vspace{-0.25cm}
\end{equation}
and Equation~\eqref{eqB2} becomes: \vspace{-0.15cm}
\begin{equation}
\label{eqdB2}
b^{(r)}_i \gets \frac{1}{d(s_i)} \left( r(s_i, \pi(s_i)) + \sum_{s_j \in S_r \setminus \{s_i\}} P^{(\pi)}(s_i, s_j) \, b^{(r)}_j \right). \vspace{-0.2cm}
\end{equation}

We now present Algorithm \ref{algo-Eval}, which provides a summary of the complete policy evaluation procedure described for both discounted reward and average reward criterion. \vspace{-0.2cm}
\begin{algorithm}[ht]
\SetAlgoLined
\SetKwInOut{Input}{Input}
\SetKwInOut{Output}{Output}
\caption{\label{algo-Eval} Policy evaluation for policy $\pi$, SISDMC-SC model}
\Input{Transition matrix $P^{(\pi)} $, reward $r(s,\pi(s)) $, partitions $\{S_r\}_{r=1}^K $, criterion: \texttt{average} or $\gamma$ \texttt{discounted}}
\Output{Value function $V^{(\pi)} $; average reward $\rho^{(\pi)} $ if \texttt{average} criterion}
\BlankLine

\If{criterion == \texttt{discounted}}{
    Update transition matrix: $P^{(\pi)} \gets \gamma \cdot P^{(\pi)} $ 
}

\If{criterion == \texttt{average}}{
    Execute \textit{\textbf{Algorithm~\ref{algo-ChiuNew}}} to compute $\rho^{(\pi)} $
}

\ForEach{partition $S_r = \{s_{1,r}, \dots, s_{n_r,r}\} $}{
    - Identify superstate $s_{1,r} $, and sets $S_{r,R} $, $S_{r,\overline{R}} $; see Equation~\eqref{eqSets} \\
    - Build local system $M^{(r)} V_{\text{sup}} = b^{(r)} $ via bottom-up recursive substitution:

    \If{criterion == \texttt{average}}{
        Use Equation \eqref{eqM1}, \eqref{eqB1}, \eqref{eqM2} and \eqref{eqB2}
    }
    \If{criterion == \texttt{discounted}}{
        Use Equation \eqref{eqM1}, \eqref{eqdB1}, \eqref{eqM2} and \eqref{eqdB2}
    }
    }
    Extract the superstates system; Equation \eqref{eq:Vsup1} and \eqref{eq:Vsup2} \\
     \If{criterion == \texttt{average}}{
     Fix reference value (e.g., $V(s_{1,1}) = 0 $) and solve the reduced linear system \\ 
     }
     \If{criterion == \texttt{discounted}}{
     Solve the reduced linear system 
     }
    \ForEach{partition $S_r = \{s_{1,r}, \dots, s_{n_r,r}\} $}{
    Reconstruct $V(s) $ for all $s_{i,r} \in S_r \setminus \{s_{1,r}\} $ using Equation \eqref{eqReconstruct}}
\end{algorithm}
\vspace{-0.2cm}
\subsection{Complexity analysis}
\begin{lemma}
\label{lemComplexite2}
    The computational complexity of the proposed policy evaluation procedure (Algorithm~\ref{algo-Eval}) is \vspace{-0.3cm}
\begin{equation}
\begin{array}{ll}
\textbf{Average reward:} & \mathcal{O}\left( \sum_{r=1}^K m_r + N \cdot K + K^3 \right) \\
\textbf{Discounted reward:} & \mathcal{O}\left( N \cdot K + K^3 \right)
\end{array}
\vspace{-0.15cm}
\end{equation}
Where $m_r $ is the number of transitions within partition $S_r $. This is significantly more scalable than classical value evaluation methods, which typically involve solving a system over all $N $ states with complexity $\mathcal{O}(N^3) $. In structured SISDMDP models where $K \ll N $, the proposed approach yields a substantial computational advantage.\vspace{-0.2cm}
\end{lemma}
\begin{proof}
Let $S_r $ denote a partition containing approximately $n_r \approx \frac{N}{K} $ states.\vspace{-0.2cm}

\paragraph{Step 0 (Average Reward Computation).}
This step is only required under the \texttt{average} reward criterion. The average reward $\rho^{(\pi)} $ is computed using Algorithm~\ref{algo-ChiuNew}, with complexity $\mathcal{O}\left( \sum_{r=1}^K m_r + K^3 \right) $, as established in Lemma~\ref{lemComplexite1}. \vspace{-0.2cm}

\paragraph{Step A (Local Substitution).}
For each partition, a local system $M^{(r)} V_{\text{sup}} = b^{(r)} $ is constructed using bottom-up recursive substitution. Due to the structured dependencies in SISDMC-SC, each local construction costs $\mathcal{O}(n_r \cdot K) $, resulting in a total cost of $\mathcal{O}(N \cdot K) $ across all partitions. This step is performed under both criteria. \vspace{-0.2cm}

\paragraph{Step B (Global System Extraction).}
One equation per partition (corresponding to the superstate) is extracted to obtain a reduced system of size $K \times K $. This operation has a cost of $\mathcal{O}(K^2) $. \vspace{-0.2cm}

\paragraph{Step C (System Resolution).}
The reduced system is solved using a direct method such as Gauss-Jordan elimination, with complexity $\mathcal{O}(K^3) $. Alternatively, iterative solvers may be used with cost $\mathcal{O}(iter\cdot K^2)$, which negligible when $K \ll N $. \vspace{-0.2cm}

\paragraph{Step D (Final Injection).}
For each non-superstate $s \in S_r \setminus \{s_{1,r}\} $, the value is reconstructed via a linear combination involving up to $K $ terms. Across all states, this step has complexity $\mathcal{O}(N \cdot K) $.
By summing all steps, we obtain the overall complexity stated above.
\end{proof}
\vspace{-0.2cm}

We now recall the Policy Iteration (PI) algorithm, integrating our structure-based policy evaluation scheme into its evaluation step. In the next section, we present numerical comparisons between this modified PI algorithm, classical PI using standard evaluation methods, and other baseline approaches such as the Value Iteration (VI).
\begin{algorithm}[ht]
\SetAlgoLined
\SetKwInOut{Input}{Input}
\SetKwInOut{Output}{Output}
 \caption{\label{algo-PI} Modified Policy Iteration, SISDMDP}
 \Input{State space $S$; action space $A$; transition matrices $P^{(a)}$; reward matrices $R^{(a)}$; criterion: \texttt{average} or $\gamma$ \texttt{discounted}}
 \Output{Optimal policy $\pi^*$, value function $V^{(\pi^*)}$}
 \BlankLine
    Set $k \gets 1$   \\
    Select an arbitrary policy $\pi_k$ \\
    \textbf{Policy evaluation}:  \\
    - Perform the SISDMC-SC evaluation procedure \textit{\textbf{Algorithm~\ref{algo-Eval}}} for $\pi_k $, according to the selected criterion \\
    - Receive $V^{(\pi_k)} $ and, if \texttt{average} criterion, $\rho^{(\pi_k)} $ \\
     \textbf{Policy Improvement}:  \\
     \If{criterion == \texttt{average}}{
     Set $\lambda \gets 1$
     }
    \If{criterion == \texttt{discounted}}{
     Set $\lambda \gets \gamma$
     }
     - Compute the Q-value from Equation \eqref{eqQ} \\
     - Choose a new policy $\pi_{k+1}$ using Equation \eqref{eqpi*} \\
     \textbf{Stopping criteria}:  \\
     \eIf{$\pi_{k+1}(s) = \pi_{k}(s)$, $\forall s \in S$}{
     -Set $\pi^*(s) \gets \pi_{k}(s)$ \\
     -Compute problem-specific average performance metrics (e.g., expected release, delay, energy ...etc) \\
     -The algorithm stops.
     }
     {
     Set $k \gets k+1$, and go to \textbf{Policy evaluation} step.
     }
\end{algorithm}
     \vspace{-0.15cm}
\begin{lemma}
\label{lemComplexite3}
    The computational complexity of the overall modified policy iteration algorithm (Algorithm~\ref{algo-PI}) is: \vspace{-0.3cm}
\begin{equation}
\begin{array}{ll}
\textbf{Average reward:} &  \mathcal{O}\Big[ \ iter_{rpi}\cdot \Big( |A|\cdot N^2 + \big( \sum_{r=1}^K m_r + N \cdot K + K^3 \big) \Big) \Big],\\[1ex]
\textbf{Discounted reward:} & \mathcal{O}\Big[ \ iter_{pi}\cdot \Big( |A|\cdot N^2 +  N \cdot K + K^3 \Big) \Big].\vspace{-0.3cm}
\end{array}
\end{equation}
\end{lemma}

\begin{proof}
    The modified policy iteration algorithm alternates between two main steps until convergence. The most computationally expensive step is the \textit{policy evaluation}, whose complexity is given in Lemma~\ref{lemComplexite2}. The second step, \textit{policy improvement}, requires $\mathcal{O}(|A| \cdot N^2)$ operations, corresponding to the maximization over actions in the Q-function \eqref{eqQ} for all states. These two steps are repeated iteratively until convergence, depending on the optimization criterion: $iter_{rpi}$ iterations for the average reward case (Relative Policy Iteration), and $iter_{pi}$ iterations for the discounted reward case (Policy Iteration).
    \vspace{-0.25cm}
\end{proof}
\begin{remark}[Semi-MDP generalization]
The SISDMDP considered in this work can be naturally extended to the semi-Markov setting. A discrete-time Semi-Markov Decision Process (SMDP) is a generalization of the standard Markov Decision Process in which actions may require a variable amount of time to complete~\cite{dietterich2000}. Under any stationary policy, the induced process preserves the same transition structure as in the SISDMDP case. As a result, the structural decomposition exploited by our procedure remains fully applicable. The only required adaptation is the inclusion of a multiplicative adjustment based on the expected holding times, for both average and discounted criteria. \vspace{-0.5cm}
\end{remark}
\section{Numerical results \label{sec:results}}
\vspace{-0.2cm}
To evaluate the performance of the proposed method, we present a numerical comparison under both criteria. The experiments are conducted on synthetic generated SISDMDPs ranging from small to large-scale instances.

For the \textbf{average reward} case (Table~\ref{tabAVG}), we compare five algorithms. The first two, \textbf{MRPI+Chiu+GTH} and \textbf{MRPI+Chiu+RB}, are variants introduced in this work (Algorithm~\ref{algo-PI}). Both rely on the Modified Relative Policy Iteration framework combined with Chiu’s decomposition. The difference lies in the linear system solvers used during policy evaluation: MRPI+Chiu+GTH employs the GTH algorithm for all systems (both intra- and inter-superstate), whereas MRPI+Chiu+RB uses the Rob-B method for intra-superstate systems and GTH only for the inter-superstate system. The remaining algorithms are \textbf{RVI} (Relative Value Iteration), \textbf{RPI+FP} (Relative Policy Iteration with Fixed-Point iterative policy evaluation), and \textbf{RPI+GJ} (Relative Policy Iteration with Gauss-Jordan elimination in the policy evaluation step).

In the \textbf{discounted reward} case (Table~\ref{tabDSC}), we compare four algorithms: \textbf{VI} (Value Iteration), \textbf{PI+FP} (Policy Iteration with Fixed-Point evaluation), \textbf{PI+GJ} (Policy Iteration with Gauss-Jordan evaluation), and our proposed method \textbf{MPI+Chiu+RB}, adapted to the discounted setting (Algorithm~\ref{algo-PI}).

Note that the difference between MRPI+Chiu+RB and MRPI+Chiu+GTH lies solely in the computation of the steady-state probability distribution. However, the computation of $V^{(\pi)}$ is identical in both cases, following the same proposed approach. This also explains the exclusion of MPI+Chiu+GTH in the discounted setting, which does not require the steady-state distribution. 

\vspace{-0.15cm}
\paragraph{Synthetic SISDMDPs generation:}
Each SISDMDP is generated from three input parameters: the total number of states $N$, the number of superstates $K$, and the action space size $|A|$. We first partition the state space into $K$ disjoint subsets of equal size $n_r=N/K$. Within each partition, one root state is designated, and a directed acyclic structure is constructed by randomly selecting forward neighbors (including the possibility of loops and revisiting previously assigned nodes). The transitions are ordered from lower to higher indexed states to ensure a hierarchical structure. Then, backward arcs are added to introduce cyclicity at the local level. To ensure connectivity at the global level, we construct a directed cycle among the $K$ superstates. Additional transitions are also introduced between states across partitions as well as among superstates themselves. While all partitions contain the same number of states, the local structure of transitions may vary due to randomly controlled transitions, resulting in diverse local dynamics. However, the randomness is controlled via consistent probabilistic rules, ensuring reproducibility for any given $(|A|, N, K)$ configuration (see source code \cite{SAYA25} for details). For instance, with $N = 10^5$ and $K = 100$, the total number of transitions across all partitions satisfies $\sum_{r=1}^K m_r \approx 719640$. Once a well-structured transition matrix is generated for the first action, the transition matrices for the remaining actions are obtained by randomly perturbing the initial probabilities, followed by normalization to preserve valid distributions.
 Instant rewards are also randomly generated for each state-action pair. As stated earlier, such structures can naturally emerge in real-world systems, particularly those governed by periodic behaviors. For instance, in \cite{YAJM24,YAJM25}, each state models a discrete number of energy packets, along with other features such as time of day or PhotoVoltaic failure status, in an energy storage system. Actions correspond to probabilistic energy transfers (e.g., selling or supplying batteries to neighboring networks). Similarly, in \cite{YHJM18,YHJM19}, states represent the number of SDUs (Service Data Units) within an optical container, following similarly structured and stochastic dynamics.

\paragraph{Stopping criteria \cite{Putr94}:}
For the average reward setting, the stopping criterion used in both RVI and in the iterative policy evaluation step of RPI+FP is based on the span seminorm, i.e., $\text{span}(V_{k+1}^{(\pi)} - V_k^{(\pi)}) < \epsilon$, or until a maximum number of iterations is reached. In contrast, for the discounted reward case, the stopping condition relies on the $\ell_\infty$ norm $\|V_{k+1}^{(\pi)} - V_k^{(\pi)} \|_\infty < \epsilon$, which leverages the contraction property of the Bellman operator under a discount factor $\gamma < 1$. We set $\epsilon = 10^{-15}$ and $\texttt{MaxIterations} = 10^5$ in all experiments. However, in large-scale instances under the average reward setting, we observed oscillations in the span value that could hinder convergence. To mitigate this, we employed a stagnation window of 100 iterations with a stagnation threshold of $10^{-13}$.
\vspace{-0.15cm}
\paragraph{Performance analysis.}
Tables~\ref{tabAVG} and~\ref{tabDSC} report the execution times (in seconds) and the number of iterations required for convergence under the average and discounted reward criteria, respectively. Each table presents two scenarios: a moderate-scale case with up to $|A| = 200$ actions and $N = 5000$ states, and a large-scale case with $|A| = 1000$ actions and up to $N = 10^5$ states. We also vary the number of partitions $K \in \{10,\ 100\}$. Each $(|A|, N, K)$ configuration is evaluated through a single run.\footnote{Execution times varied by less than ±10\% over 30 randomized runs for a fixed configuration, based on 95\% confidence intervals.} The fastest algorithm for each configuration is also highlighted.

\vspace{-0.1cm}
\begin{itemize}
    \item In both average and discounted reward settings, all algorithms based on policy iteration or relative policy iteration (RPI+FP, RPI+GJ, MRPI+Chiu+RB, etc.) require the same number of iterations for a given configuration. The advantage of our methods lies in accelerating the \textit{policy evaluation} step, which dominates the computational cost. For example, in the average reward case (Table~\ref{tabAVG}), with $|A| = 1000$, $N = 10^5$, and $K = 10$, MRPI+Chiu+RB converges in 1105.75 seconds using an exact solver, compared to $2899.64$ seconds for RPI+FP (a fixed-point method that may be less precise), despite both requiring $7$ iterations. Other methods exceed $10^4$ seconds in this configuration. Similarly, in the discounted setting (Table~\ref{tabDSC}), PI+Chiu+RB solves the largest instance in 233.12 seconds, whereas PI+FP takes 2101.06 seconds. It is also worth noting that Rob-B-based approaches are even faster in the discounted setting, mainly because they avoid computing the steady-state probability distribution.

    \item The impact of $K$ is more significant in our decomposable methods (MRPI+Chiu+RB, MRPI+Chiu+GTH and PI+Chiu+RB), where $K$ explicitly appears in the complexity expressions (Lemma \ref{lemComplexite3}). A larger $K$ reduces the size of each partition ($N/K$ states), which limits the benefits of our propagation mechanism. Conversely, smaller values of $K$ lead to larger partitions, which can still be handled efficiently by our method. For instance, with $|A| = 1000$ and $N = 10^5$, MRPI+Chiu+RB takes 532.80 seconds for $K = 100$ (5 iterations) versus 233.12 seconds for $K = 10$ (6 iterations). This supports our design assumption that the efficiency of our method improves when $K \ll N$.

    \item Value iteration methods (RVI and VI) remain competitive in moderate-scale scenarios (top sections of Tables~\ref{tabAVG} and~\ref{tabDSC}), particularly when the number of partitions is high ($K = 100$). This is due to the internal propagation overhead of our method, which increases as $K$ grows. However, value iteration struggles to scale in larger instances, given its overall complexity of $O(\text{iter}_{vi} \cdot |A| \cdot N^2)$.

    \item The RPI+GJ and PI+GJ approaches are clearly limited to moderate-scale problems, as solving the linear system in the evaluation step has cubic complexity. The same limitation applies to MRPI+Chiu+GTH, which uses the GTH algorithm in all subsystems. This becomes especially problematic when $N$ is large and $K$ is small, making each subsystem (of size $N/K$) expensive to solve. For example, for $K = 10$ and $N = 10^4$, MRPI+Chiu+GTH requires 5989.38 seconds, and for larger systems, execution time exceeds $10^4$ seconds.
\end{itemize}
Overall, these results strongly support the effectiveness of the proposed methods, MRPI+Chiu+RB and PI+Chiu+RB, which consistently deliver exact solutions with substantial runtime improvements in large-scale SISDMDPs, especially when $K\ll N$.
\paragraph{Source code:}
\vspace{-0.2cm}
Algorithms were implemented using a Python-based framework specifically developed for this work \cite{SAYA25}, with efficient handling of sparse matrices via vectorized operations. Experiments were conducted on a laptop equipped with 10 CPU cores (8 cores at 3.2~GHz peak frequency and 2 cores at 2.0~GHz), and 16~GB of RAM.  
\begin{table}[ht]
\centering
\caption{\textbf{\textit{Average reward criterion}} – Exec. time (s) and number of iterations.\label{tabAVG}}
\vspace{0.2cm}

%%%%%%%%%%%%%%%%%%%%%%%%%%%%%%%%%%%%%%%%%%%%%%%%%%%%
% ------------------ A = 200 ---------------------- %
%%%%%%%%%%%%%%%%%%%%%%%%%%%%%%%%%%%%%%%%%%%%%%%%%%%%
\begin{tabular}{|c|c||cc|cc|cc}
\multicolumn{2}{||c||}{} & \multicolumn{6}{c}{$|A| = 200$} \\
\cmidrule{3-8}
 & \textbf{Algorithm} 
 & \multicolumn{2}{c|}{$N=10^{3}$} 
 & \multicolumn{2}{c|}{$N=3\!\times\!10^{3}$} 
 & \multicolumn{2}{c}{$N=5\!\times\!10^{3}$} \\
\cline{3-8}
 &  
 & $K=100$ & $K=10$ 
 & $K=100$ & $K=10$ 
 & $K=100$ & $K=10$ \\
\midrule
\multirow{2}{*}{}
 & RVI            & \cellcolor{yellow!30}0.68    & 2.66     & \cellcolor{yellow!30}4.54    & 3.58     & \cellcolor{yellow!30}10.82    & 5.77     \\
 &                & \cellcolor{yellow!30}201     & 1013     & \cellcolor{yellow!30}492     & 643      & \cellcolor{yellow!30}751      & 818      \\
\multirow{2}{*}{}
 & RPI+FP         & 4.25     & 10.67   & 20.30    & 25.67   & 39.88    & 25.21    \\
 &                & 6        & 5       & 6        & 5       & 6        & 5        \\
\multirow{2}{*}{}
 & RPI+GJ         & 12.30    & 10.17   & 141.53   & 119.61  & 492.65   & 325        \\
 &                & 6        & 5       & 6        & 5       & 6        & 5        \\
 \multirow{2}{*}{}
 & MRPI+Chiu+GTH      & 8.47     & 7.62    & 29.61    & 175.97  & 61.59    & 763.89   \\
 &                & 6        & 5       & 6        & 5       & 6        & 5        \\
\multirow{2}{*}{}
 & MRPI+Chiu+RB   & 8.25     & \cellcolor{yellow!30}0.89    & 16.83    & \cellcolor{yellow!30}2.81    & 25.80    & \cellcolor{yellow!30}5.20     \\
 &                & 6        & \cellcolor{yellow!30}5       & 6        & \cellcolor{yellow!30}5       & 6        & \cellcolor{yellow!30}5        \\
\bottomrule
\end{tabular}

\vspace{0.25cm}
%\addlinespace[0.2cm]

%%%%%%%%%%%%%%%%%%%%%%%%%%%%%%%%%%%%%%%%%%%%%%%%%%%%
% ------------------ A = 1000 --------------------- %
%%%%%%%%%%%%%%%%%%%%%%%%%%%%%%%%%%%%%%%%%%%%%%%%%%%%
\begin{tabular}{|c|c||cc|cc|cc}
\multicolumn{2}{||c||}{} & \multicolumn{6}{c}{$|A| = 1000$} \\
\cmidrule{3-8}
 & \textbf{Algorithm} 
 & \multicolumn{2}{c|}{$N=10^{4}$} 
 & \multicolumn{2}{c|}{$N=5\!\times\!10^{4}$} 
 & \multicolumn{2}{c}{$N=10^{5}$} \\
\cline{3-8}
 &  
 & $K=100$ & $K=10$ 
 & $K=100$ & $K=10$ 
 & $K=100$ & $K=10$ \\
\midrule
\multirow{2}{*}{}
 & RVI            & 183.09  & 68.67   & 1507.79   & 846.92  & $>10^4$ & $>10^4$  \\
 &                & 1185   & 476      & 1333      & 744     & 817     & 664      \\
\multirow{2}{*}{}
 & RPI+FP         & 156.34    & 95.16  & 645.51   & 656.48   & 2858.74  & 2899.64  \\
 &                & 6         & 6     & 5        & 8         & 7        & 7       \\
\multirow{2}{*}{}
 & RPI+GJ         & 2505.44  & 2717.38 & $>10^4$  &$>10^4$   & $>10^4$ & $>10^4$ \\
 &                & 6        & 6       & 5       & 8        & 7        &  7       \\
 \multirow{2}{*}{}
 & MRPI+Chiu+GTH  & 162.24    & 5989.38 &$>10^4$ &$>10^4$    & $>10^4$ & $>10^4$  \\
 &                & 6         & 6     & 5        & 8         & 7       & 7        \\
\multirow{2}{*}{}
 & MRPI+Chiu+RB   & \cellcolor{yellow!30}40.75   & \cellcolor{yellow!30}12.06  & \cellcolor{yellow!30}258.41   & \cellcolor{yellow!30}234.93  & \cellcolor{yellow!30}1368.82 & \cellcolor{yellow!30}1105.75   \\
 &                & \cellcolor{yellow!30}6       & \cellcolor{yellow!30}6      & \cellcolor{yellow!30}5        & \cellcolor{yellow!30}8       & \cellcolor{yellow!30}7        & \cellcolor{yellow!30}7        \\
\bottomrule
\end{tabular}
\end{table}
\begin{table}[ht]
\centering
\caption{\textbf{\textit{Discounted ($\gamma = 0.9$) reward criterion}} – Exec. time (s) and number of iterations. \label{tabDSC}}
\vspace{0.2cm}

%%%%%%%%%%%%%%%%%%%%%%%%%%%%%%%%%%%%%%%%%%%%%%%%%%%%
% ------------------ A = 200 ---------------------- %
%%%%%%%%%%%%%%%%%%%%%%%%%%%%%%%%%%%%%%%%%%%%%%%%%%%%
\begin{tabular}{|c|c||cc|cc|cc}

\multicolumn{2}{||c||}{} & \multicolumn{6}{c}{$|A| = 200$} \\
\cmidrule{3-8}
 & \textbf{Algorithm} 
 & \multicolumn{2}{c|}{$N=10^{3}$} 
 & \multicolumn{2}{c|}{$N=3\!\times\!10^{3}$} 
 & \multicolumn{2}{c}{$N=5\!\times\!10^{3}$} \\
\cline{3-8}
 &  
 & $K=100$ & $K=10$ 
 & $K=100$ & $K=10$ 
 & $K=100$ & $K=10$ \\
\midrule
\multirow{2}{*}{}
 & VI            & \cellcolor{yellow!30}1.34   & 1.13   & \cellcolor{yellow!30}2.65   & 1.71   & \cellcolor{yellow!30}4.83     & 2.90    \\
 &                & \cellcolor{yellow!30}365    & 341   & \cellcolor{yellow!30}371    & 340    & \cellcolor{yellow!30}362      & 359     \\
\multirow{2}{*}{}
 & PI+FP         & 4.73   & 4.39   & 11.31   & 10.08  & 22.41  & 16.94   \\
 &                & 6      & 6      & 5       & 5      & 6      & 5      \\
\multirow{2}{*}{}
 & PI+GJ         & 10.81    & 11.51  & 105.84  & 106.11 & 431.28 & 366.46  \\
 &                & 6       & 6     & 5       & 5      & 6      & 5       \\
\multirow{2}{*}{}
 & MPI+Chiu+RB   & 1.78   & \cellcolor{yellow!30}0.31   & 4.02    & \cellcolor{yellow!30}0.71   & 8.01   & \cellcolor{yellow!30}1.24    \\
 &                & 6      & \cellcolor{yellow!30}6      & 5       & \cellcolor{yellow!30}5      & 6      & \cellcolor{yellow!30}5       \\
\bottomrule
\end{tabular}

\vspace{0.25cm}
%\addlinespace[0.2cm]

%%%%%%%%%%%%%%%%%%%%%%%%%%%%%%%%%%%%%%%%%%%%%%%%%%%%
% ------------------ A = 1000 --------------------- %
%%%%%%%%%%%%%%%%%%%%%%%%%%%%%%%%%%%%%%%%%%%%%%%%%%%%
\begin{tabular}{|c|c||cc|cc|cc}

\multicolumn{2}{||c||}{} & \multicolumn{6}{c}{$|A| = 1000$} \\
\cmidrule{3-8}
 & \textbf{Algorithm} 
 & \multicolumn{2}{c|}{$N=10^{4}$} 
 & \multicolumn{2}{c|}{$N=5\!\times\!10^{4}$} 
 & \multicolumn{2}{c}{$N=10^{5}$} \\
\cline{3-8}
 &  
 & $K=100$ & $K=10$ 
 & $K=100$ & $K=10$ 
 & $K=100$ & $K=10$ \\
\midrule
\multirow{2}{*}{}
 & VI            & 147.17     & 55.63 s     & 417.50   & 374.41     & $>10^4$    & $>10^4$   \\
 &                & 354       & 357         & 364      & 355        & 374        & 366       \\
\multirow{2}{*}{}
 & PI+FP         & 57.32    & 42.02     & 180.49     & 213.58     & 1379.09     & 2101.06   \\
 &                & 5        & 5         & 5          & 6          & 5          & 6         \\
\multirow{2}{*}{}
 & PI+GJ         & 3873.59  & 3847.86   & $>10^4$    & $>10^4$    & $>10^4$    & $>10^4$   \\
 &                & 5        & 5         & 5          & 6          & 5          & 6         \\
\multirow{2}{*}{}
 & MPI+Chiu+RB   & \cellcolor{yellow!30}20.68 & \cellcolor{yellow!30}5.21 &  \cellcolor{yellow!30} 73.90 &  \cellcolor{yellow!30} 26.43 &  \cellcolor{yellow!30} 532.80 &  \cellcolor{yellow!30} 233.12 \\
 &                &  \cellcolor{yellow!30} 5         &  \cellcolor{yellow!30}5         &  \cellcolor{yellow!30}5          &  \cellcolor{yellow!30}6          &  \cellcolor{yellow!30}5          &  \cellcolor{yellow!30}6         \\
\bottomrule
\end{tabular}

\end{table}
\vspace{-0.3cm}
\section{Conclusion \label{sec:conclusion}}
\vspace{-0.3cm}
In this work, we introduced the SISDMDP framework, a structured class of Markov Decision Processes that leverages single-input decompositions and recurrence properties to enable efficient policy evaluation. Building on this structure, we proposed exact solution methods applicable to both average and discounted reward settings. The proposed algorithms significantly reduce computation time in large-scale MDPs while maintaining full accuracy, particularly by accelerating the policy evaluation step.
Our numerical experiments demonstrate the scalability and effectiveness of the approach. Beyond algorithmic contributions, the SISDMDP model offers a promising direction for the structured modeling of real-world decision systems, such as multi-station battery management or queueing systems with spatial partitioning. An interesting direction for future work is to explore how this structure can be incorporated into model-free reinforcement learning. In particular, integrating SISDMDP-compatible decompositions into Q-learning or deep RL frameworks could enable more efficient learning in large and structured environments.
\vspace{-0.4cm}
\section*{Acknowledgment}
\vspace{-0.2cm}
This work is partially supported by the public grant of the Fondation Math\'ematique Jacques Hadamard (FMJH) through the PGMO-UVSQ program.

\bibliographystyle{splncs04}
\bibliography{MDP_structured}

\begin{thebibliography}{10}
\providecommand{\url}[1]{\texttt{#1}}
\providecommand{\urlprefix}{URL }
\providecommand{\doi}[1]{https://doi.org/#1}

\bibitem{YHJM18}
{Ait EL Mahjoub}, Y., {Castel-Taleb}, H., Fourneau, J.M.: Performance and
  energy efficiency analysis in ngreen optical network. In: 14th International
  Conference on Wireless and Mobile Computing, Networking and Communications
  (WiMob) (2018). \doi{10.1109/WiMOB.2018.8589144}

\bibitem{YHJM19}
{Ait EL Mahjoub}, Y., {Castel-Taleb}, H., Fourneau, J.M.: A numerical approach
  of the analysis of optical container filling. In: 12th EAI ValueTools (2019).
  \doi{10.1145/3306309.3306333}

\bibitem{YAJM24}
{Ait El Mahjoub}, Y., Fourneau, J.M.: Finding the optimal policy to provide
  energy for an off-grid telecommunication operator. In: 20th International
  Conference on Wireless and Mobile Computing, Networking and Communications
  (WiMob) (2024). \doi{10.1109/WiMob61911.2024.10770514}

\bibitem{YAJM25}
{Ait El Mahjoub}, Y., Fourneau, J.M.: A slot-based energy storage
  decision-making approach for optimal off-grid telecommunication operator.
  Computer Communications journal  (2025). \doi{10.1016/j.comcom.2025.108273}

\bibitem{SAYA25}
Alouah, S., {Ait El Mahjoub}, Y.: {SISDMDP Framework} - source code (2025),
  \url{https://github.com/ossef/SISDMDP_Research}

\bibitem{barto2003}
Barto, A.G., Mahadevan, S.: Recent advances in hierarchical reinforcement
  learning. Discrete Event Dynamic Systems  \textbf{13}(4) (2003)

\bibitem{boutilier1995}
Boutilier, C., Dearden, R., Goldszmidt, M.: Exploiting structure in policy
  construction. In: Proc. 14th International Joint Conference on Artificial
  Intelligence (IJCAI) (1995),
  \url{https://www.ijcai.org/Proceedings/95-2/Papers/012.pdf}

\bibitem{buchholz1994}
Buchholz, P.: Exact and ordinary lumpability in finite markov chains. Journal
  of Applied Probability  \textbf{31}(1) (1994). \doi{10.2307/3215235}

\bibitem{courtois1977}
Courtois, P.J.: Decomposability: Queueing and Computer System Applications.
  Academic Press (1977)

\bibitem{dietterich2000}
Dietterich, T.G.: Hierarchical reinforcement learning with the maxq value
  function decomposition. Journal of Artificial Intelligence Research
  \textbf{13} (2000). \doi{doi.org/10.1613/jair.639}

\bibitem{Chiu87}
Feinberg, B.N., Chiu, S.S.: A method to calculate steady-state distributions of
  large markov chains by aggregating states. Operations Research  (1987).
  \doi{10.1287/opre.35.2.282}

\bibitem{Muntz1994}
Franceschinis, G., Muntz, R.R.: Bounds for quasi-lumpable markov chains.
  Performance Evaluation  \textbf{20} (1994).
  \doi{10.1016/0166-5316(94)90015-9}, performance '93

\bibitem{Abhj15}
Gosavi, A.: Simulation-Based Optimization: Parametric Optimization Techniques
  and Reinforcement Learning. Springer New York, NY (2015).
  \doi{doi.org/10.1007/978-1-4899-7491-4}

\bibitem{GTHe85}
Grassman, W., Taksar, M., Heyman, D.: Regenerative analysis and steady state
  distributions for {Markov} chains. Operations Research  \textbf{33}(5),
  1107--1116 (1985)

\bibitem{koller1999}
Koller, D., Parr, R.: Computing factored value functions for policies in
  structured mdps. In: Proc. 16th International Joint Conference on Artificial
  Intelligence (IJCAI) (1999). \doi{10.5555/646307.687921}

\bibitem{Marin2022}
Marin, A., Piazza, C., Rossi, S.: Proportional lumpability and proportional
  bisimilarity. Acta Informatica  \textbf{59} (2022).
  \doi{10.1007/s00236-021-00404-y}

\bibitem{Putr94}
Puterman, M.L.: Markov Decision Processes: Discrete Stochastic Dynamic
  Programming. John Wiley \& Sons, Inc. (1994). \doi{10.1002/9780470316887}

\bibitem{Rob90}
Robertazzi, T.G.: Computer Networks and Systems: Queueing Theory and
  Performance Evaluation. Springer New York, NY (1990).
  \doi{doi.org/10.1007/978-1-4684-0385-5}

\bibitem{SLTD17}
Song, Y., Lin, J., Tang, M., Dong, S.: An internet of energy things based on
  wireless lpwan. Engineering  \textbf{3}(4) (2017).
  \doi{10.1016/J.ENG.2017.04.011}

\bibitem{Stew94}
Stewart, W.J.: Introduction to the Numerical Solution of Markov Chains.
  Princeton University Press (1994)

\bibitem{sutton1999}
Sutton, R.S., Precup, D., Singh, S.: Between mdps and semi-mdps: A framework
  for temporal abstraction in reinforcement learning. Artificial Intelligence
  \textbf{112} (1999)

\bibitem{LAZM15}
Vangelista, L., Zanella, A., Zorzi, M.: Long-range iot technologies: The dawn
  of lora™. In: Future Access Enablers of Ubiquitous and Intelligent
  Infrastructures (09 2015). \doi{10.1007/978-3-319-27072-2\_7}

\end{thebibliography}

\end{document}